\documentclass{gtmon_a}
\pdfoutput=1
\input{diagxy}
\let\to\rightarrow
\let\barrsquare\square
\let\square\undefined

\usepackage[T1,T5]{fontenc}

\proceedingstitle{Proceedings of the School and Conference in Algebraic
Topology (The Vietnam National University, Hanoi, 9-20 August 2004)}
\conferencestart{9 August 2004}
\conferenceend{20 August 2004}
\conferencename{School and Conference in Algebraic Topology}
\conferencelocation{Vietnam National University, Hanoi, Vietnam}

\editor{John Hubbuck}
\givenname{John}
\surname{Hubbuck}

\editor{Nguy\~\ecircumflex{}n H\,V H\uhorn{}ng}
\givenname{Nguy\~\ecircumflex{}n H V}
\surname{Hu'ng}

\editor{Lionel Schwartz}
\givenname{Lionel}
\surname{Schwartz}

\title{Algebraic Topology (Hanoi, August 2004)\\The problem session}

\author[Broto, H\uhorn{}ng, Kuhn, Palmieri, Priddy and Yagita]{Carles
    Broto, {\fontencoding{T5}\selectfont
    Nguy\~\ecircumflex{}n H\,V H\uhorn{}ng}, Nicholas J
    Kuhn\newline John H Palmieri, Stewart Priddy, Nobuaki Yagita}
\givenname{Carles}
\surname{Broto}
\address{Departament de Matem\`atiques\\
Universitat Aut\`onoma de Barcelona\\\newline
08193 Bellaterra\\
Spain}
\email{broto@mat.uab.es}
\urladdr{}

\givenname{Nguy\~\ecircumflex{}n H V}
\surname{Hu'ng}
\address{Department of Mathematics\\
Vietnam National University\\\newline
334 Nguyen Trai Street\\
Hanoi\\
Vietnam}
\email{nhvhung@vnu.edu.vn}
\urladdr{}

\givenname{Nicholas J}
\surname{Kuhn}
\address{Department of Mathematics\\
University of Virginia\\\newline
Charlottesville VA 22903\\
USA}
\email{njk4x@virginia.edu}
\urladdr{}

\givenname{John H}
\surname{Palmieri}
\address{Department of Mathematics\\
University of Washington\\\newline
Box 354350\\
Seattle WA 98195-4350\\
USA}
\email{palmieri@math.washington.edu}
\urladdr{}

\givenname{Stewart}
\surname{Priddy}
\address{Department of Mathematics\\
Northwestern University\\\newline
Evanston IL 60208\\
USA}
\email{priddy@math.northwestern.edu}
\urladdr{}

\givenname{Nobuaki}
\surname{Yagita}
\address{Department of Mathematics\\
Faculty of Education\\
Ibaraki University\\\newline
Mito\\
Ibaraki\\
Japan}
\email{yagita@mx.ibaraki.ac.jp}
\urladdr{}

\volumenumber{11}
\issuenumber{}
\publicationyear{2007}
\papernumber{20}
\startpage{435}
\endpage{441}

\doi{}
\MR{}
\Zbl{}

\arxivreference{}


\received{24 November 2006}
\revised{}
\accepted{}
\published{14 November 2007}
\publishedonline{14 November 2007}
\proposed{}
\seconded{}
\corresponding{}
\version{}

\makeatletter
\def\cnewtheorem#1[#2]#3{\newtheorem{#1}{#3}[section]
\expandafter\let\csname c@#1\endcsname\c@Theorem}

\let\xysavmatrix\xymatrix
\def\xymatrix{\disablesubscriptcorrection\xysavmatrix}
\AtBeginDocument{\let\bar\wbar\let\tilde\wtilde\let\hat\what}
\AtBeginDocument{}

\newtheorem{Theorem}{Theorem}[section]
\cnewtheorem{Proposition}[Theorem]{Proposition}
\cnewtheorem{Lemma}[Theorem]{Lemma}
\cnewtheorem{Corollary}[Theorem]{Corollary}
\cnewtheorem{Conjecture}[Theorem]{Conjecture}
\theoremstyle{remark}
\cnewtheorem{Definition}[Theorem]{Definition}
\cnewtheorem{Remark}[Theorem]{Remark}
\cnewtheorem{Question}[Theorem]{Question}
\cnewtheorem{Observation}[Theorem]{Observation}
\cnewtheorem{Problem}[Theorem]{Problem}
\cnewtheorem{Example}[Theorem]{Example}
\cnewtheorem{Subsec}[Theorem]{}
\cnewtheorem{Note}[Theorem]{}

\makeatother  

\makeautorefname{Conjecture}{Conjecture}

\def\eex{{\accent"5E e}\kern-.470em\raise.3ex\hbox{\char'176}}
\def\uw{u\kern-.44em\raise.82ex\hbox{ \vrule width .12em height .0ex 
depth .075ex \kern-0.16em \char'56}\kern-.05em}
\def\EEX{{\accent"5E E}\kern-.60em\raise.9ex\hbox{\char'176}
\kern.1em}
\def\UW{U\kern-.42em\raise1.36ex\hbox{
\vrule width .13em height .0ex depth .075ex \kern-0.16em
\char'56}\kern-.07em}
\def\aah{{\accent"5E a}\kern-.62em\raise.2ex\hbox{\char'22}\kern.12em}

\newcommand{\V}{{\mathbb{V}}}
\newcommand{\F}{{\mathbb{F}}} 
\newcommand{\gls}{\mathrm{GL}_s}
\newcommand{\otigls}{\otimes_{\gls}}
\newcommand{\N}{\mathbb{N}}
\makeop{hocolim}
\newcommand{\Ext}{\rm{Ext}}
\newcommand{\Fp}{\mathbb{F}_{p}}
\newcommand{\ext}[1]{\Ext_{P}^{#1}(\Fp,\Fp)}
\newcommand{\p}{\mathcal{P}^{0}}

\makeop{im}
\makeop{res}
\makeop{det}
\makeop{Tr}
\makeop{tr}
\makeop{Sq}
\makeop{colim}

\begin{document}

\begin{abstract}
This article contains a collection of problems contributed during the
course of the conference.
\end{abstract}

\maketitle

\section{The problems presented by Carles Broto}
Broto, Levi and Oliver \cite{BLO2,BLO-surv} introduced the concept of
$p$--local finite group as an algebraic object modeled on the fusion on the
Sylow $p$--subgroup of a finite group, and attached to it a classifying
space which is a $p$--complete space with properties similar to that of
the $p$--completed classifying space of a finite group. Every finite group
$G$ gives rise canonically to a $p$--local finite group with classifying
space homotopy equivalent to $BG^\wedge_p$.  But there are also exotic
examples; that is, $p$--local finite groups that cannot be obtained from
a finite group in the canonic way and therefore its classifying space
is not homotopy equivalent to the $p$--completed classifying space of
any finite group.

Some exotic examples are obtained by Broto and M{\o}ller \cite{BM}
out of $p$--compact groups (see Dwyer--Wilkerson \cite{DW-pc} and M{\o}ller
\cite{Moeller}). More precisely, to a 1--connected
$p$--compact group $X$ it is attached a family $X(q)$
of $p$--local finite groups, $q$ a $p$--adic unit,   that approximates $X$ in the sense that the classifying spaces
$\{ BX(q^{p^m})\}_{m\in \N}$ form a direct system with mapping telescope
$$\hocolim_{m}BX(q^{p^m}) \simeq BX$$
(up to $p$--completion). In particular, this shows that every 1--connected
$p$--compact group can be expressed, up to $p$--completion, as a mapping
telescope of $p$--local finite groups. 
In view of that, Clarence Wilkerson asked the following  question:

\begin{Question}
Is every $p$--compact group a telescope of finite groups? 
In more precise words, given a $p$--compact group $X$, does there exist
a sequence of finite groups and homomorphisms
$$\bigl\{ G_0\stackrel{{}_{f_1}}{\longrightarrow}\cdots
  \stackrel{{}_{f_i}}{\longrightarrow} G_i
  \stackrel{{}_{f_{i+1}}}{\longrightarrow}\cdots \bigr\}$$
and a homotopy equivalence $\hocolim_{i}BG_i \simeq BX$, up to $p$--completion?
\end{Question}

A possibly related question is the following.

\begin{Question}
Let  $\{ BX_m\}_{m\in \N}$ be a direct system
of classifying spaces of $p$--local finite groups and continuous maps 
$f_m\colon BX_m\to BX_{m+1}$ with 
$\hocolim_{m}BX_m \simeq BX$, the classifying space of a $p$--compact group $X$. 
Is there a map $\hat g\colon BG^\wedge_p\to BX_m $, for some $m$, making the diagram 
$$\bfig
  \morphism|b|<400,0>[\cdots`BX_m;f_m]
  \morphism(400,0)|b|<400,0>[\phantom{BX_m}`\cdots;f_{m+1}]
  \morphism(800,0)<350,0>[\phantom{\cdots}`BX;]
  \morphism(1150,400)|r|<0,-400>[BG^{\wedge}_p`\phantom{BX};g]
  \morphism(1150,400)<-750,-400>[\phantom{[BG^{\wedge}_p}`\phantom{BX_m};\hat{g}]
  \efig$$
homotopy commutative, 
where $BG^\wedge_p$ is the $p$--completed classifying space of a finite group
(or $p$--local finite group)?
\end{Question}


\section{The problem presented by Nguy\~\ecircumflex{}n H\,V H\uhorn{}ng}
Let $\Sq^0\co \text{Ext}^{s,t}_{\mathcal{A}}(\F_2, \F_2) 
\to \text{Ext}^{s,2t}_{\mathcal{A}}(\F_2, \F_2)$
be the squaring operation induced on the cohomology of 
the Steenrod algebra by the Frobenius map 
$\mathcal{A}_* \to \mathcal{A}_*, x \mapsto x^2$.

\begin{Conjecture}
\label{Sq^0onExt}
($\Sq^0$ is eventually isomorphic on the Ext groups.)

Let $\text{Im}(\Sq^0)^i$ denote the image of $(\Sq^0)^i$ on
$\text{Ext}_{\mathcal{A}}^s(\F_2, \F_2)$.  Then, for arbitrary $s$,
there exists a number $r=r(s)$ depending on $s$ such that
$$(\Sq^0)^{i-r}\co \text{Im}(\Sq^0)^r \to \text{Im}(\Sq^0)^i$$
  is an isomorphism for every $i \geq r$. 
\end{Conjecture}

In other words, the conjecture predicts that any finite $\Sq^0$--family
in $\text{Ext}_{\mathcal{A}}^s(\F_2, \F_2)$ has at most $r(s)$ nonzero
elements.

\begin{Observation}
We guess that $r(s) = s-2$. 
\end{Observation}

Let us explain the motivation for this conjecture.

Denote by $\V_s$ an $s$--dimensional vector space over $\F_2$, and let 
$$\Tr_s\co \F_2\otigls PH_d(B\V_s)\to \text{Ext}_{\mathcal{A}}^{s,s+d}(\F_2,\F_2)$$
be the algebraic transfer defined by W~Singer as an
algebraic version of the geometrical transfer 
$\tr_s\co \pi_*^S((B\V_s)_+) \to \pi_*^S(S^0)$. Here $PH_*(B\V_s)$ denotes the
primitive part of $H_*(B\V_s)$ consisting of all homology classes that are
annihilated by any Steenrod operations of positive degrees.

It has been proved that $\Tr_s$ is an isomorphism for $s=1,2$ by
   Singer and for $s=3$ by Boardman.
   Among other things, these data together with the fact that
   $\Tr=\oplus_{s} \Tr_s$ is an algebra homomorphism 
show that $\Tr_s$ is highly nontrivial. Therefore,
   the algebraic transfer is expected to be a useful tool in the
   study of the mysterious cohomology of the Steenrod algebra,
   $\text{Ext}_{\mathcal{A}}^{*,*}(\F_2,\F_2)$.
Although Singer recognizes that $\Tr_5$ is not isomorphic, 
his open conjecture predicts that $\Tr_s$ is a monomorphism for any $s$. 

According to Boardman and Minami, one has a commutative diagram
$$\bfig
  \barrsquare<1500,500>[(\F_2{\otigls}PH_*(B\V_s))_d`
    \mathrm{Ext}_{\mathcal{A}}^{s,s+d}(\F_2,\F_2)`
    (\F_2{\otigls}PH_*(B\V_s))_{2d+s}`
    \mathrm{Ext}_{\mathcal{A}}^{s,2(s+d)}(\F_2,\F_2)~,;
    \Tr_s`\Sq^0`\Sq^0`\Tr_s]
  \efig$$
where the left vertical arrow is the Kameko $\Sq^0$ and the right
vertical one is the classical squaring operation.

The conjecture comes from the above diagram and the following theorem.

\begin{Theorem}[{{\fontencoding{T5}\selectfont H\uhorn{}ng \cite{Hung}}}]
\label{stable}
($\Sq^0$ is eventually isomorphic on the domain of the transfer.)

Let $\text{Im}(\Sq^0)^i$ denote the image of the Kameko iterated 
squaring $(\Sq^0)^i$ on the domain of the transfer, $\F_2\otigls PH_*(B\V_s)$.
Then 
$$(\Sq^0)^{i-s+2}\co \text{Im}(\Sq^0)^{s-2} \to \text{Im}(\Sq^0)^i$$
is an isomorphism for every $i\geq s-2$.
\end{Theorem}

\section{The problems presented by Nick Kuhn}

These three problems are related to pondering the Singer Transfer
$$\tau_{s,t}\co \operatorname{Hom}_{{\mathcal  A}}^{t-s}(H^*(({\mathbb Z}/p)^s;{\mathbb F}_p),{\mathbb F}_p)_{GL_s({\mathbb Z}/p)} \rightarrow \operatorname{Ext}_{{\mathcal  A}}^{s,t}({\mathbb F}_p,{\mathbb F}_p)$$

\begin{Problem}
Prove or disprove the conjecture that $\operatorname{Ext}_{{\mathcal  A}}^{s,t}({\mathbb F}_2,{\mathbb F}_2)= 0$ if $\alpha(t)>s$, where, as usual, $\alpha(t)=k$ if $t = 2^{i_1}+ \dots + 2^{i_k}$ with $i_1< \dots i_k$.
\end{Problem}

Note that Reg Wood's verification of the Peterson conjecture shows that the domain of $\tau_{s,t}$ is 0 if $\alpha(t)>s$, and I am conjecturing that the range of $\tau_{s,t}$ vanishes in these same places. Bob Bruner has checked that the conjecture is true in the range of known $\operatorname{Ext}$ calculations. \\

\begin{Problem}
Can $\tau_{*,*}$ be viewed as the edge homomorphism of some spectral sequence?

This problem is not so well posed, but the point would be to find some
sensible way to measure the failure of $\tau_{*,*}$ to be an isomorphism
of bigraded algebras.
\end{Problem}

\begin{Problem}
Fix $s$.  Is there a uniform calculation of $\operatorname{Ext}_{{\mathcal  A}}^{s,*}({\mathbb F}_p,{\mathbb F}_p)$ for all large enough $p$?
\end{Problem}

For example, $\operatorname{Ext}_{{\mathcal  A}}^{1,*}({\mathbb F}_p,{\mathbb F}_p)$ can be viewed as `independent of $p$' for $p>2$. The heuristic here goes as follows: $\operatorname{Ext}_{{\mathcal  A}}^{s,*+s}({\mathbb F}_p,{\mathbb F}_p)$ is similar to the representation theoretic object $\operatorname{Hom}_{{\mathcal  A}}^{*}(H^*(({\mathbb Z}/p)^s;{\mathbb F}_p),{\mathbb F}_p)_{GL_s({\mathbb Z}/p)}$, and thus is perhaps `independent of $p$' for large $p$ analogous to similar situations studied in the work of Anderson, Jantzen, and Soergel.


\section{The problems presented by John H Palmieri}

Let $p$ be a prime.  Let $P_{*}\cong \Fp [\xi_{1}, \xi_{2}, \xi_{3},
\dots]$ be the polynomial part of the dual Steenrod algebra, and write
$P$ for its dual; thus when $p=2$, $P$ is the mod 2 Steenrod algebra,
and when $p$ is odd, $P$ is a quotient Hopf algebra of the mod $p$
Steenrod algebra.  Consider $\ext{*,*}$.

The $p$th power map on $P_{*}$ induces an algebra endomorphism $\p$ on
Ext, which with respect to the grading acts like this:
\[
\p \co \ext{s,t} \rightarrow \ext{s,pt}.
\]

\begin{Conjecture}\label{conj}
Fix $s>0$.  For all $z \in (\p)^{-1} \ext{s,t}$, $z$ is nilpotent.
Equivalently, for all $z \in \ext{s,t}$, there exists an $n\geq 0$ so
that $(\p)^{n} z$ is nilpotent.
\end{Conjecture}

The following is proved in \cite{palmieri;quotient}.

\begin{Theorem}
\fullref{conj} is true when $p=2$.
\end{Theorem}

Hand and computer calculations when $p=2$ suggest the following.

\begin{Conjecture}
In \fullref{conj}, one may take $n$ to be 1.  That is, for all
$z \in \ext{s,t}$, $\p z$ is nilpotent.
\end{Conjecture}

\section{The problems presented by Stewart Priddy}

Let $G_n = GL_n(\mathbb{F}_p)$. Quillen has shown 
$\colim_n H_i(G_n;\mathbb{F}_p)=0$ for $i>0$. The classes of 
$H^*(G_n;\mathbb{F}_p)$ are unstable characteristic classes for 
representations over $\mathbb{F}_p$. 

Now let $p=2$; Maazen has shown $H^i(G_{2n};\mathbb{F}_2) = 0$ for $0 <i<n$.
What about $i=n$? A maximal elementary $p$--subgroup of $G_{2n}$ is
$$ A_{n,n}=
\left(\begin{matrix} I_n & *   \\ 0 & I_n \end{matrix} \right).$$
$A_{n,n} \approx (\mathbb{Z}/2)^{n^2}$ and 
$$H^*(A_{n,n}; \mathbb{F}_2) = \mathbb{F}[x_{i,j} | 1 \leq i,j \leq n ]$$
where $x_{i,j}$ are the obvious 1--dimensional generators. 
The normalizer
$N_{G_{2n}}(A_{n,n}) = G_n \times G_n$; thus the restriction homomorphism
has the form
$$\res\co H^*(G_{2n}; \mathbb{F}_2) \longrightarrow 
H^*(A_{n,n};\mathbb{F}_2)^{G_n \times G_n}$$
We consider the image in dimension $n$. Let 
$$\det_n = \det{x_{i,j}}$$
Then $\det_n \in H^n(A_{n,n};\mathbb{F}_2)^{G_n \times G_n}$

\begin{Question}
Is $\det_n \in \im(\res)$?
\end{Question}

This is true for $n=1,2$ (see Milgram and Priddy \cite{mp}). 
There is also an analogous question for $p$ odd but one must
take into account the action of the diagonal matrices.

\begin{Question}
Is $H^n(G_{2n} ; \mathbb{F}_2) = \mathbb{F}_2$
generated by $\det_n$?
\end{Question}

\section{The problem presented by Nobuaki Yagita}
Let us write by $E$ the extraspecial $p$--group $p_+^{1+2}$
of order $p^3$ and exponent $p$ for an odd prime $p$.
In Theorem 6.2 of my paper in the present proceedings 
(Stable splitting and cohomology of
$p$--local finite groups over the extraspecial $p$--group
of order $p^3$ and exponent $p$), we have a graph which shows
stable splitting of $BG$ for all (p-local) finite groups having
a Sylow $p$--subgroup $E$ for $p=3$. The problem is to write down the
similar graph for $p=7$.  A partial result is given in Theorem 9.4
in the same paper.

\bibliographystyle{gtart}
\bibliography{link}

\end{document}